\documentclass[12pt]{amsart}
\usepackage{amsmath,amssymb,graphicx,subfigure,psfrag}
\newtheorem{thm}{Theorem}
\newtheorem{lem}[thm]{Lemma}

\newtheorem{alg}[thm]{Algorithm}
\newtheorem{rules}[thm]{Rules}

\newcommand{\bi}{\begin{itemize}}
\newcommand{\ei}{\end{itemize}}
\newcommand{\be}{\begin{enumerate}}
\newcommand{\ee}{\end{enumerate}}
\newcommand{\bc}{\begin{center}}
\newcommand{\ec}{\end{center}}
\newcommand{\bt}{\begin{tabular}}
\newcommand{\et}{\end{tabular}}
\newcommand{\ba}{\begin{array}}
\newcommand{\ea}{\end{array}}

\newcommand{\noi}{\noindent}

\begin{document}

\title[Permutations generated by stacks in series]
{Permutations generated by a stack of\\
depth 2 and an infinite stack in series}

\author[Murray Elder] {Murray Elder} 
\address{Department of Mathematics, Stevens Institute of Technology,  Hoboken NJ USA}
\email[url]{murrayelder@gmail.com, http://www.math.stevens.edu/$\sim$melder}

 \date{\today}

\begin{abstract}
We prove that the set of permutations generated by a stack of depth two and an infinite stack in series
has a basis (defining set of forbidden patterns) consisting of 20 permutations of length 5, 6, 7 and 8. We prove this via a ``canonical'' generating algorithm. 
\end{abstract}

\maketitle

\section{Introduction}

\begin{figure}[ht!]
\bc
\bt{|c|c|c|}
\hline
\includegraphics[height=34.5mm]{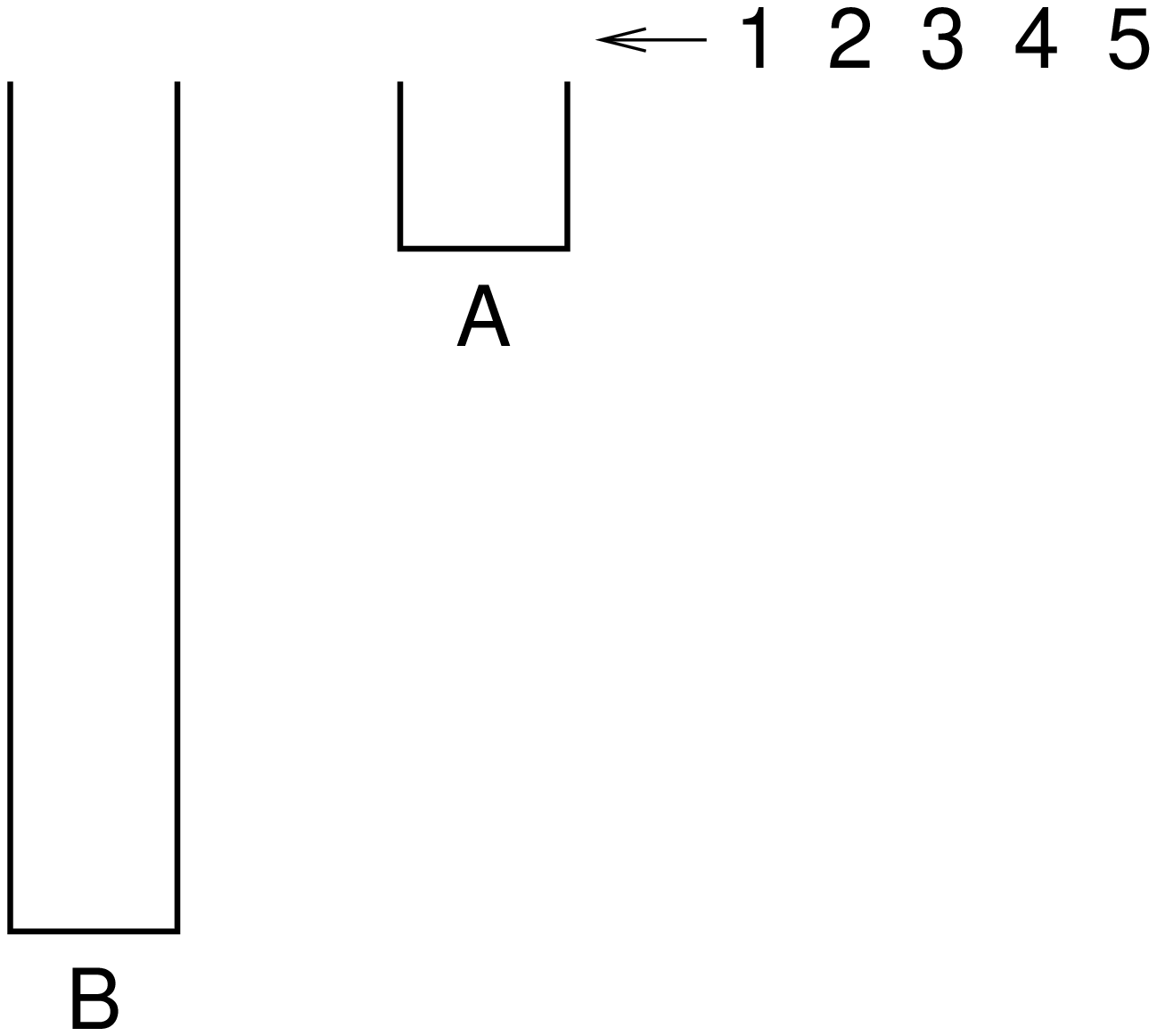} &
\includegraphics[height=34.5mm]{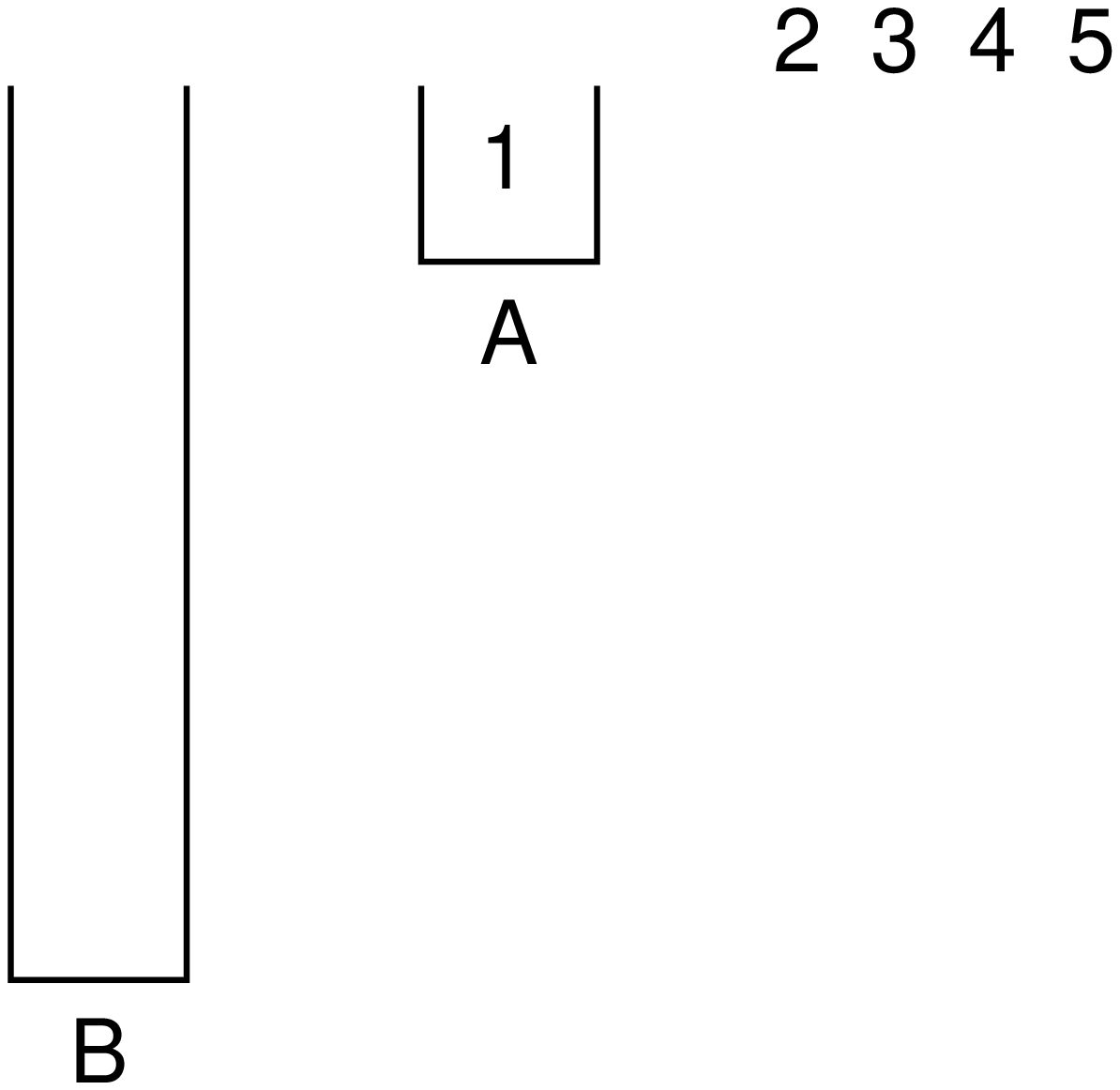} &
\includegraphics[height=34.5mm]{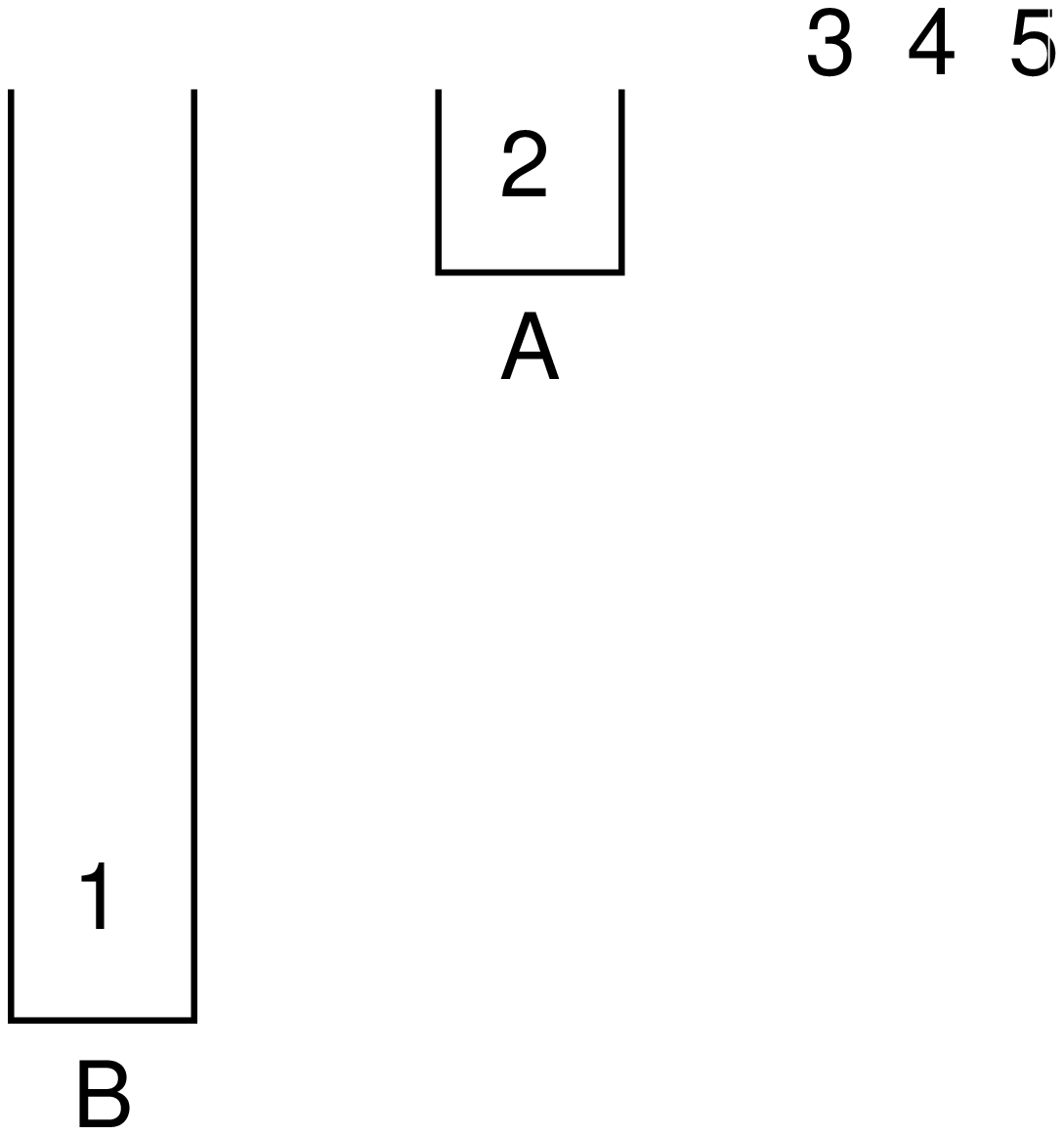} \\
\hline
& \texttt{Apply rule 1.2} & \texttt{Apply rule 3.1} \\
\hline
\includegraphics[height=34.5mm]{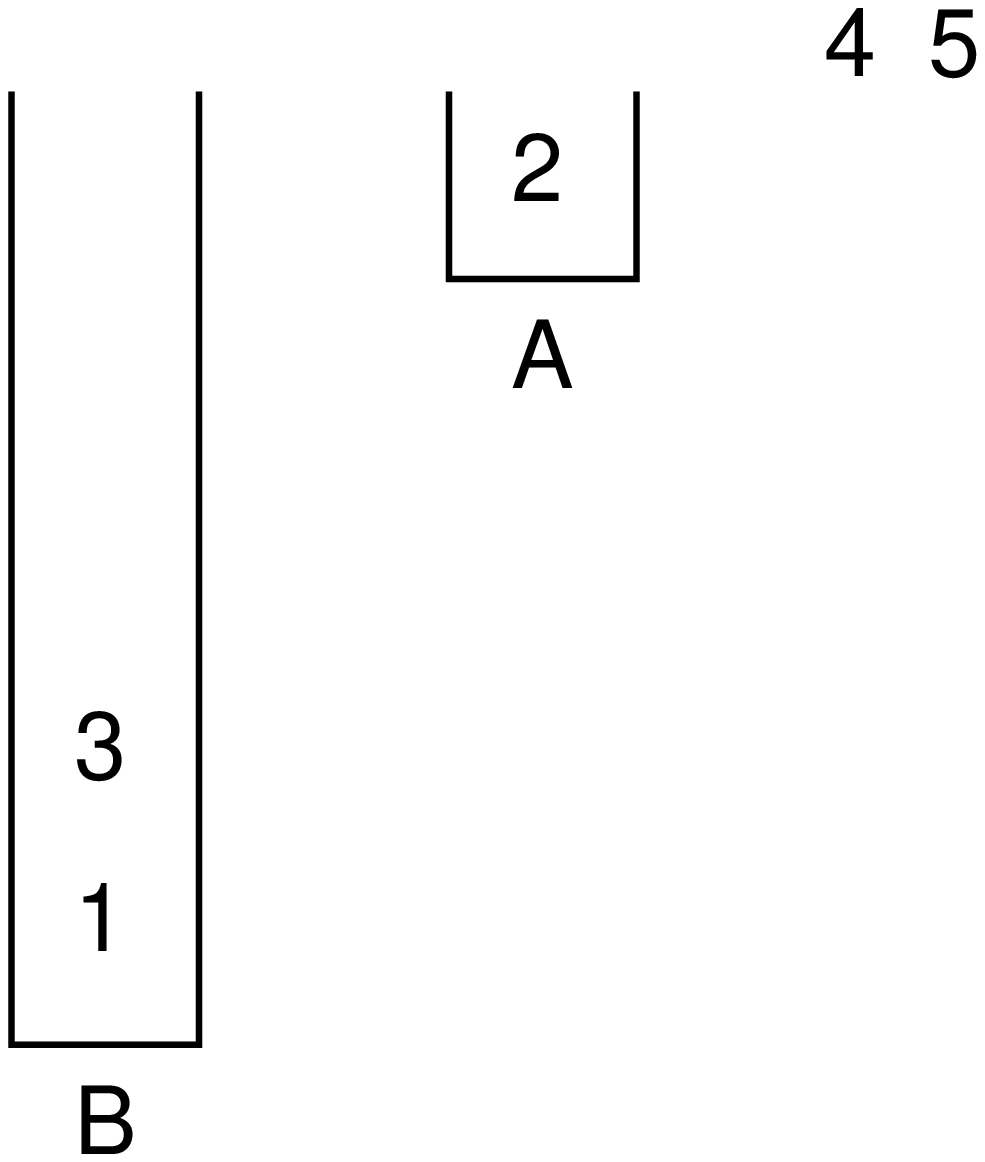}  &
\includegraphics[height=34.5mm]{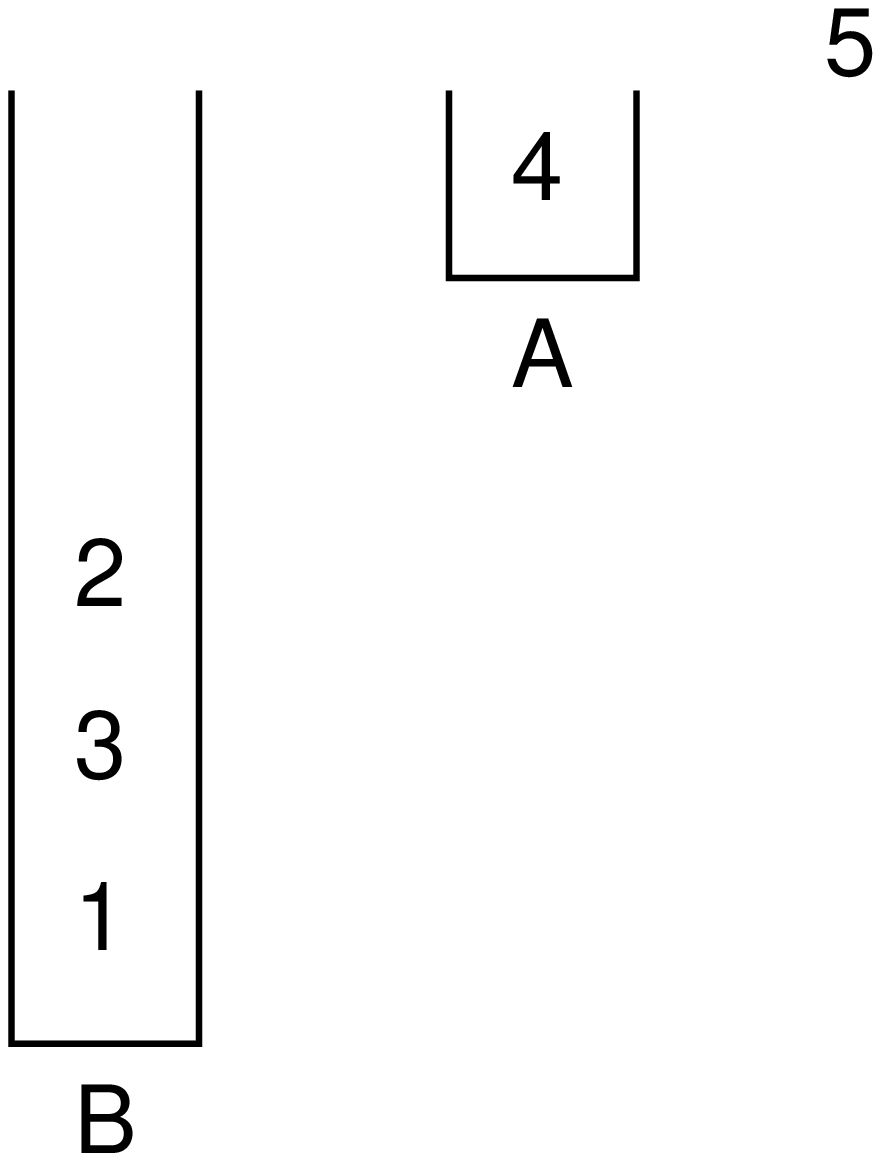} &
\includegraphics[height=34.5mm]{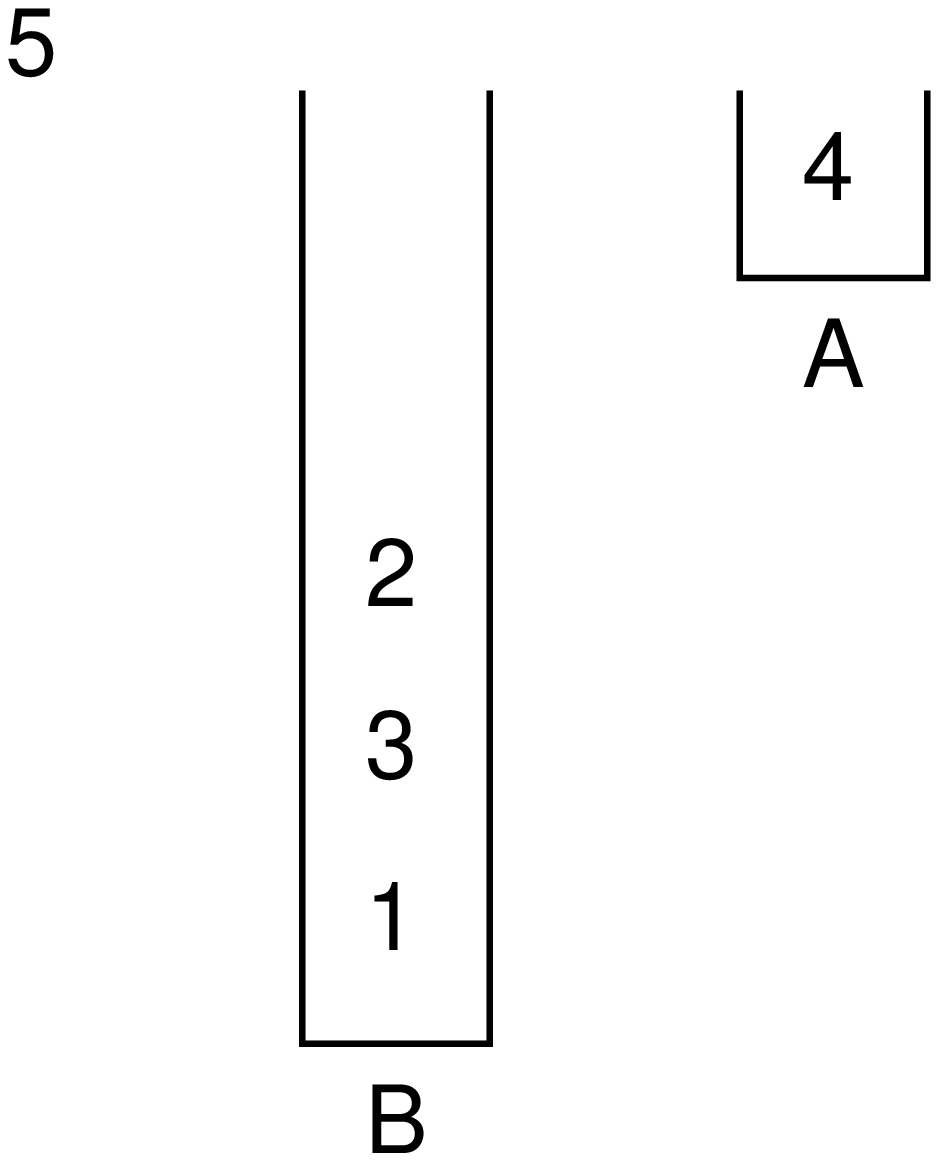} \\
\hline
\texttt{Apply rule 2.2} & & \\
\hline
\includegraphics[height=34.5mm]{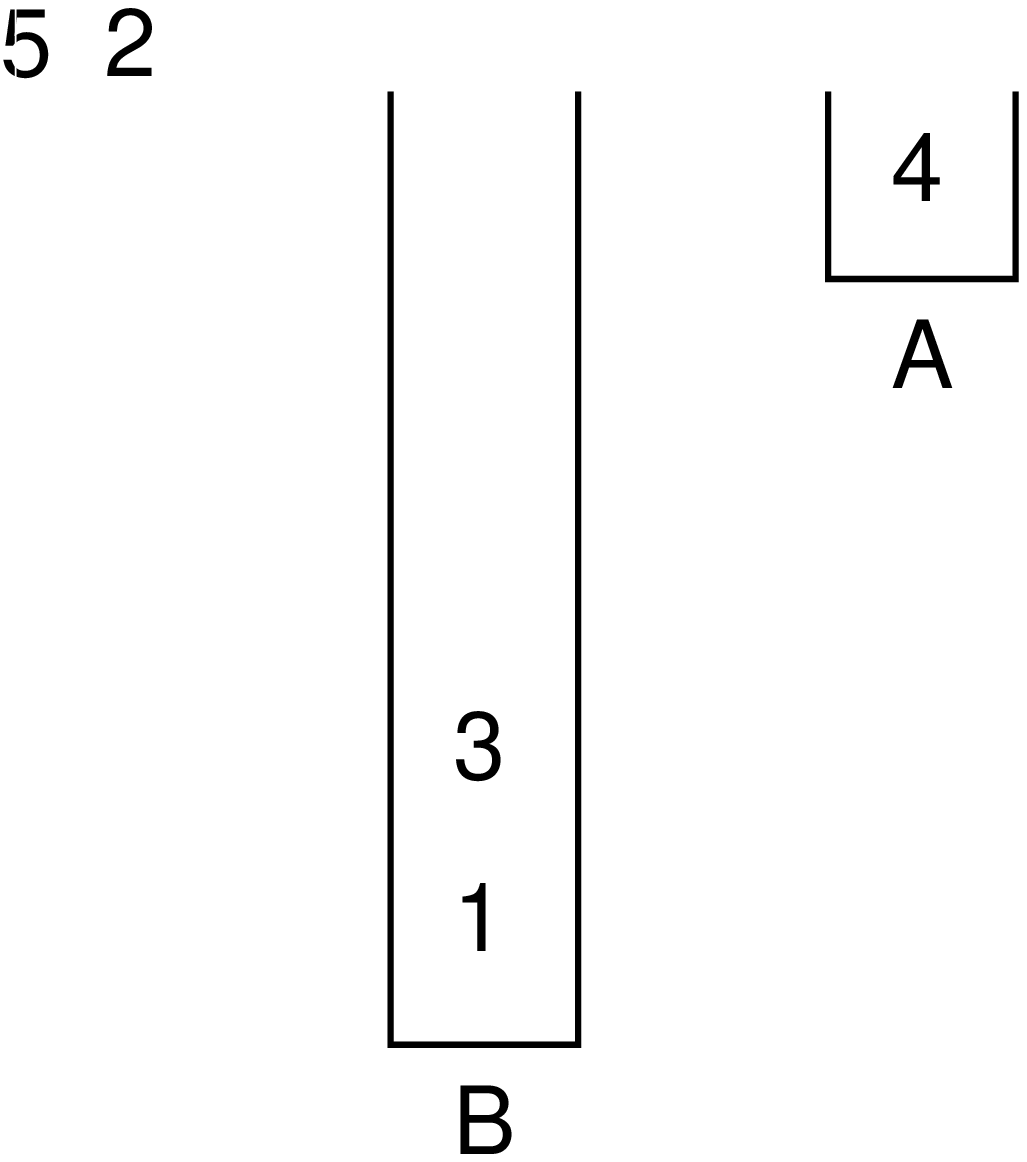} &
\includegraphics[height=34.5mm]{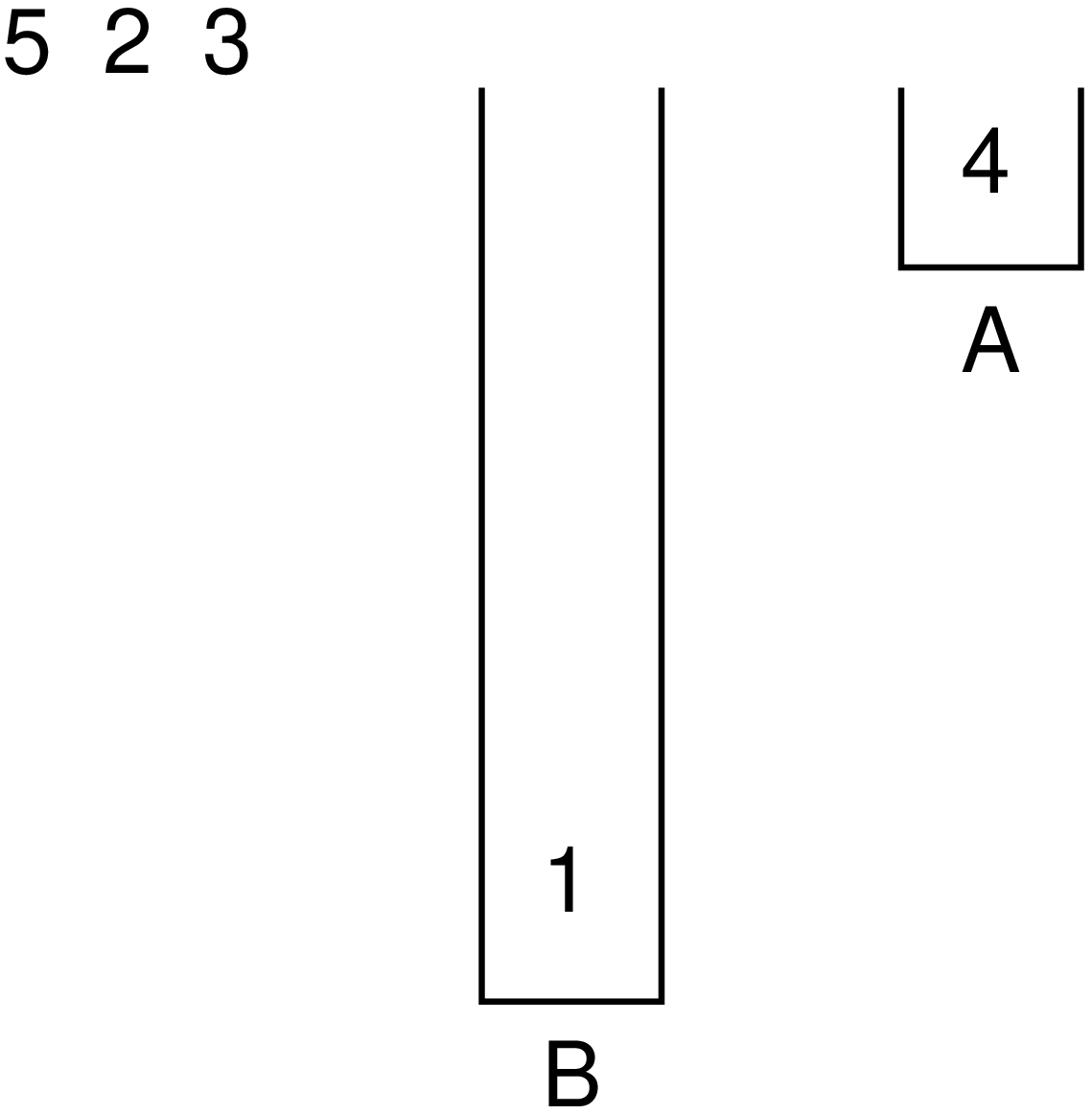} &
\includegraphics[height=34.5mm]{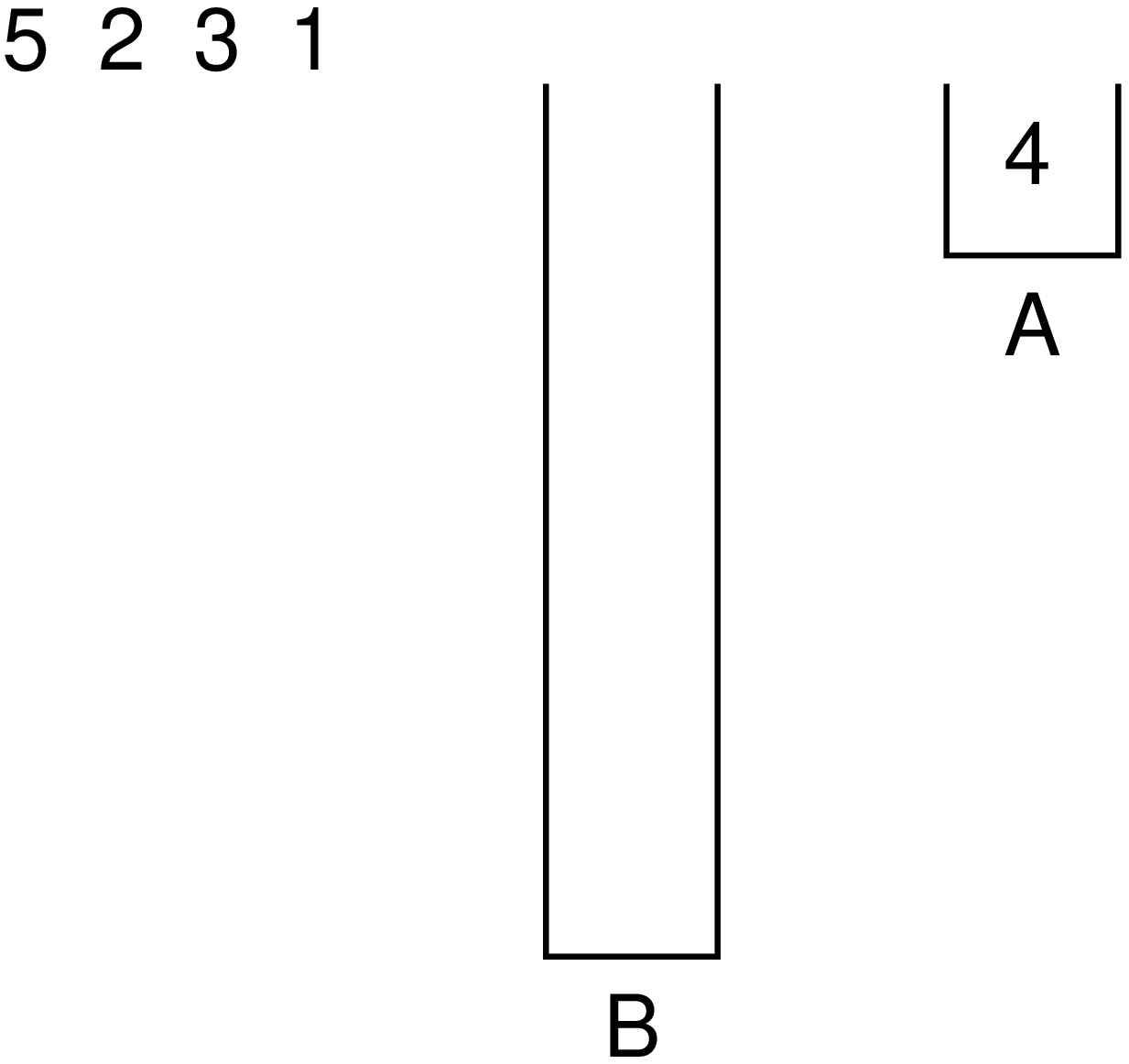} \\
\hline
\et \ec
\caption{Generating the permutation $52314$.}\label{fig:stacks}
\end{figure}

In this article we examine the set of permutations that can be generated by passing the sequence $1, 2, \ldots, n$ through a stack of depth two followed by an infinite stack, as in Figure \ref{fig:stacks}. The {\em depth} of a stack is the number of tokens it can hold, including one space at the top for passing tokens through the stack.    By convention we pass tokens right to left.
 
We prove in Theorem \ref{thm:main} that a permutation can be generated in this way if and only if it avoids a list of sub-patterns of 20 permutations, and furnish a deterministic procedure (Algorithm \ref{alg}) for generating them. These permutations were found initially by computations with a stack of depth two and a stack of depth $k$ for increasing $k$ by Linton \cite{sal}.

The current interest in permutations that avoid sub-patterns could perhaps be traced back to  Knuth, who proved that a permutation can be generated by passing an ordered sequence though a single infinite stack if and only if it avoids the subsequence $3, 1, 2$ \cite{Knuth} \footnote{Actually he proved the equivalent fact that a permutation can be sorted if and only if it avoids $2, 3, 1$. He also showed that they are enumerated by the Catalan numbers.}.

For two infinite stacks in series, the set of avoided minimal sub-patterns that characterize  the permutations that can be generated is infinite \cite{Max}. But somewhere between a first stack of depth one (\textit{ie.} no first stack) and infinite depth, there is a break point where the basis goes from being finite to infinite (see Lemma \ref{lem:basis1}).

A good overview of permutations generated and sorted in various ways using stacks can be read in 
\cite{bona-book} and a good introduction to the field of pattern avoiding permutations can be found in 
\cite{bona-stacks}. Recent open problems in the field are summarized in \cite{me-vince}.

The article is organized as follows. In Section \ref{sec:prelim} we define permutations and pattern avoidance, and give some basic facts and terminology for permutations generated by stacks. In Section \ref{sec:alg} we describe an algorithm to decide whether or not a given permutation can be generated using a stack of depth two followed by an infinite stack. We prove that the algorithm is valid, and that  a permutation  is accepted if and only if it can be generated by the stacks, if and only if it avoids the 20 permutations.

\section{Preliminaries}\label{sec:prelim}

A {\em permutation} is an arrangement of a finite number of distinct elements of a linear order, for example, $5, 1, 2, 4, 3$ or $4, 6, 1$. It is customary to omit the commas and write $51243$. Two permutations are {\em order isomorphic} if they have the same relative ordering. So $231$ and $461$ are order isomorphic. Define a {\em sub-permutation} of a permutation $p_1\ldots p_n$ to be a word $p_{i_1}\ldots p_{i_s}$
with $i_1<\ldots <i_s$. A {\em subinterval} of a permutation is a sub-permutation consisting of contiguous entries, that is, $i_{j+1}=i_j+1$ for $j=1,\ldots, s-1$.  
A permutation $p$ {\em contains} or {involves} a permutation $q$ if it has a sub-permutation that  is order isomorphic to $q$. So $p=51243$ contains $q=321$ since deleting the entries $1$ and $2$ of $p$ gives the sub-permutation $543$ which is order isomorphic to $q$. 
A permutation $p$ {\em avoids} $q$ if it does not contain it. So $51243$ avoids the permutation $231$ since  no sub-permutation is order isomorphic to $231$.

A set of permutations $S$ is said to be {\em closed} (under involvement) if $p\in S$ and $p$ involves or contains $q$ implies that $q\in S$.  Given a set of permutations $B$, the set $Av(B)$ of permutations which do not contain any permutations from $B$ is closed, and is called the {\em avoidance set} for $B$. If a set of permutations can be described as the avoidance set for some set $B$, and $B$ is the minimal such set (so that no element of $B$ contains another) then we call $B$ the {\em basis} for the set. For example, the set of permutations that avoid $12$ and $123$ is the set of all decreasing permutations, and its basis is simply $\{12\}$. Note also that if $\sigma$ is in a basis for a set $S$ then deleting any entry of $\sigma$ gives a permutation that is order isomorphic to an element of $S$.

Define $S_{k,\infty}$ to be  the set of permutations that can be generated by passing $1\ldots n$ through a stack of depth $k$ followed by an infinite stack, and define $B_{k,\infty}$ as its basis. So for example, the basis for a stack of depth one (so no storage) followed by an infinite stack is $\{312\}$.
We will call the input symbols {\em letters} or {\em tokens}.

\begin{lem}\label{lem:basis1}
Let $\sigma \in B_{k,\infty}$ and define $\sigma_3$ to be the string of integers obtained by adding three to the value of every entry of $\sigma$.
Either  $\sigma$ or  $\sigma_3 213$  is in $B_{k+1,\infty}$.
\end{lem}
\textit{Proof:}
If $\sigma$ is not in $S_{k+1,\infty}$ then, since deleting 
  any entry of $\sigma$ gives a permutation that is in 
  $S_{k,\infty}\subseteq S_{k+1,\infty}$, it follows that 
  $\sigma\in B_{k+1,\infty}$. So we can assume that 
  $\sigma$ is in $S_{k+1,\infty}$. Now consider the 
  problem of generating $\sigma_3 213$.

Whichever way you put the tokens $1, 2$ and $3$ onto the two stacks, some token must occupy the first stack. These tokens must stay until the rest of the permutation has been output, so the remaining tokens must be processed with the first stack of depth $k$ rather than $k+1$. But since $\sigma$ cannot be generated with the first stack of depth $k$ then neither can $\sigma_3$, so $\sigma_3 213$ is not in 
$S_{k+1,\infty}$.

To show that $\sigma_3 213$ is a basis element, we must show that every shorter permutation contained in it is in $S_{k+1,\infty}$. Let $\tau$ be a sub-permutation of $\sigma_3 213$ obtained by deleting one entry.  
If $\tau=\sigma_3 21,\sigma_3 23$ or $\sigma_3 13$ then we can generate it as follows. Place the first two entries $(1,2),(1,3)$ or $(2,3)$  on the second stack in the appropriate order. This leaves the first stack clear, so we can now generate $\sigma_3$ using the two stacks (since $\sigma \in S_{k+1,\infty}$), and lastly output the two tokens.  If  instead $\tau$ has an entry deleted from the $\sigma_3$ prefix, then we place tokens $1,2$ on the second stack with $2$ on top, and leave $3$ on the first stack. Since $\sigma$ was in the basis for  $S_{k,\infty}$, we can generate $\sigma_3$ with one entry deleted while the first stack has depth $k$, and then we can output $213$. Thus we can generate any sub-permutation of $\sigma_3 213$.
\hfill$\Box$

It follows that for all $n\in\mathbb N, |B_{k,\infty}|\leq  |B_{k+1,\infty}|$. Since by Theorem \ref{thm:main} $B_{2,\infty}$ is finite, then either there is a number $n> 2$ such that $B_{k,\infty}$ is finite for all $k\leq n$ and not finite for $k>n$, or $B_{n,\infty}$ could be finite for all $n\in\mathbb N$.

Define $\mathcal B$ to be  the following set of 20 permutations in Table \ref{table:B}.
\begin{table}
\[
\begin{array}{|l|r|r||l|r|r|}
 \hline
\mathrm{length} \; 5 & 51234 & 52134  & \mathrm{length} \; 7 & 4175623& 4275613\\

& 51243 & 52143 & & 4137256 & 4237156\\

 & 51423  & 52413 &  & 4137265 & 4237165\\
 \hline
\mathrm{length} \; 6 &  645123 & 645213 & \mathrm{length} \; 8    & 41386725 & 42386715\\ 

 &416235 &  426135 &&& \\
& 416253 & 426153 &&& \\
 \hline
\end{array}
\]
\caption{The set $\mathcal B$.}\label{table:B}
\end{table}
Observe that   $\mathcal B$ is closed under the operation of interchanging the  1 and 2 entries. 
\begin{lem}\label{lem:basis}
If a permutation contains an element of  $\mathcal B$ then it cannot be generated by a stack of depth $2$ followed by an infinite stack.
\end{lem}
\textit{Proof:} It suffices to prove that none of the permutations in $\mathcal B$ can be generated by the two stacks. It then follows that no permutation containing one can be generated. 
It is routine to check by hand or computer that each of the permutations in  $\mathcal B$ cannot be generated by a stack of depth two followed by an infinite stack. We can enumerate the full list of permutations of length up to 8 generated by considering codewords on three letters $\rho,\lambda,\mu$ that correspond to pushing tokens from input to the first stack $(\rho)$, from the first to the second stack $(\lambda)$, then output $(\mu)$.  For example, the codeword $\rho \lambda \mu \rho \lambda \mu \rho \lambda \mu$ generates the permutation $123$. We require that each prefix must have no more $\lambda$s than $\rho$s, no more $\mu$s than $\lambda$s, and further that the number of $\rho$s is no more than 2 plus the number of $\lambda$s.
Using this technique we can verify that none of the permutations in  $\mathcal B$ are produced.
\hfill$\Box$.

Linton conjectured that this set should be the basis for $S_{2,\infty}$. In \cite{Atk} Aktinson {\em et al}
consider the set of permutations generated by passing $1\ldots n$ through a finite  token passing network, which is a directed graph where nodes can hold at most one 
token, and tokens move in any way from an input to an output node. One can view  two  stacks in series as a token passing network, which is finite if both stacks are of bounded depth. They prove that the set of permutations generated via finite networks can be encoded in a regular language, and from this one can find its basis. Using an implementation of this procedure in GAP Linton \cite{sal} computed the bases for $S_{2,k}$ for increasing $k$, and observed that in each case 20 small permutations occurred, as well as longer permutations which related to the bound on the second stack. 
Linton conjectured that in the limit the basis should consist of just these elements.

In proving this conjecture we will make use of the following technical definitions.
Let us say that a subinterval $\tau$ of a permutation $\sigma=\alpha\tau\beta$ is {\em right-contiguous} 
if $\beta$ does not contain any entries between the minimum and maximum entries in $\tau$, and is 
{\em right-contiguous modulo $a$} if $\beta$ does not contain any entries between the minimum and maximum entries in $\tau$ except the entry $a$. For example, the subinterval $413$ of $4137256$
is right contiguous modulo $2$ and $137$ is not.

Lastly, we will make use of the following notation for permutations below. If a permutation contains a token  $a$ preceding a token $b$, then we write $-a-b-$, or simply $a-b$, when we do not know the other letters of the permutation. The notation $a<^b_c$ means that we know the permutation contains an $a$ preceding both $b$ and $c$, but  we do not know the relative orderings of the $b$ and $c$. That is, the permutation could be $a-b-c$ or $a-c-b$. For example, every permutation of length $5$ in Table \ref{table:B} except $52413$ is of the form $5<^{1-4}_2$.

\section{The algorithm}\label{sec:alg}

In this section we describe a deterministic procedure to generate permutations using a depth two followed by an infinite stack in series.
Denote the stack of depth two stack as $A$ and the infinite stack as $B$.  Without loss of generality 
if there are input letters remaining then $A$ contains an entry $a$. That is, if $A$ ever becomes empty then we will immediately fill it with the next input letter.

Let $x$ be the next input letter, and $b$ be the top entry of $B$ (if non-empty). 
The stack $B$ is {\em well  ordered} if its elements are stacked from top to bottom in order they appear 
in the permutation generated. 

\begin{figure}[ht!]
\bc
\includegraphics[height=4cm]{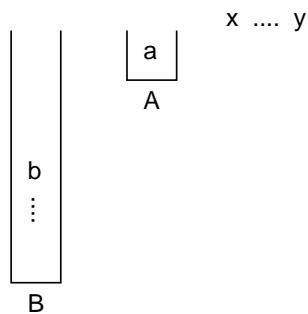}
\ec
\caption{Algorithm}\label{fig:2stacks}
\end{figure}

If $b,a$ or $x$ are next to be output, then output, and move the next input letter into $A$ if empty.
If the next output letter is below $b$ then the permutation cannot be generated (by the algorithm being used). 
Otherwise the next output letter is $y$ somewhere back in the input list. See Figure \ref{fig:2stacks}.

We will choose between pushing  $a$ to stack $B$ and $x$ to stack $A$, or vice versa, using the following (fairly technical) rules.

\begin{rules}\label{rules}
If $y$ is the next token to be output, $x$ precedes it in the input stream, $a$ is on stack $A$ and $b$ is at the top of stack $B$ (if non-empty), then:
\be
\item[1.1.] If in the output permutation $a$ precedes some two letters in the input that lie between 
$x$ and $y$, then keep $a$ on stack  $A$ and put $x$  on stack $B$.

\item[1.2.] If  in the output permutation $x$ precedes some two letters in the input that lie between 
$x$ and $y$, then put $x$ on $A$ and $a$  on $B$.
\item[2.1.] If  in the output permutation $b$ precedes $a$ then keep $a$  on $A$ and put $x$  on $B$.
\item[2.2.] If  in the output permutation $b$ precedes $x$ then put $x$  on $A$ and  $a$  on $B$.
\ee
When none of the conditions for 1.1-2.2 are met, then:
\be
\item[3.1.] If  in the output permutation $a$ precedes $x$ then keep $a$ on $A$ and put $x$ on $B$ 
{\em unless} the subinterval $y\ldots a$  is  right-contiguous modulo $x$ and avoids $312$,
 in which case put $a$ on $B$ and  $x$ on $A$.
\item[3.2.] If  in the output permutation $x$ precedes $a$ then put $a$ on $B$ and  $x$ on $A$ 
{\em unless} the subinterval $y\ldots x$ is  right-contiguous and avoids $312$, in which case keep $a$ on $A$ and  put $x$ on $B$.
\ee
\end{rules}

\begin{figure}[ht!]
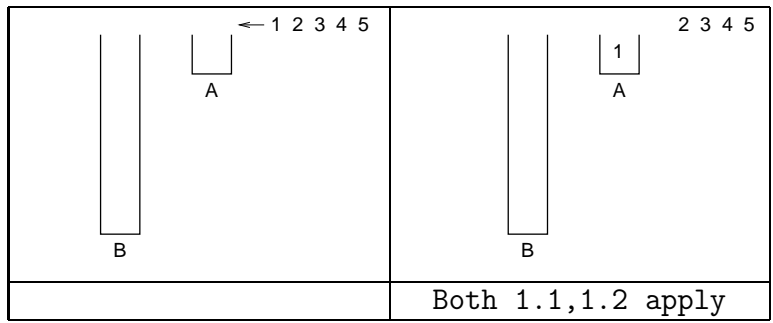

\bc
\bt{|c|c|}
\hline
\includegraphics[height= 35.0mm]{figures/fig1_1.eps} &
\includegraphics[height= 35.0mm]{figures/fig1_2.eps} \\
\hline
 &   \texttt{Both 1.1,1.2 apply}  \\
\hline
\et \ec
\caption{Generating the permutation $51234$.}\label{fig:51234}
\end{figure}

If the conditions of a rule are met, we say that the rule {\em applies}. It may happen that at some point while pushing tokens around the stacks,  two or more rules apply  simultaneously. For example, to generate the permutation $51234$, we would put $1$ on stack $A$, then $a=1,x=2$ and $y=5$. In this case, both rules 1.1 and 1.2 apply.  See Figure  \ref{fig:51234}.
If this is the case then any algorithm we construct using these rules will halt.
We now state the algorithm.

\begin{alg}\label{alg}
 Given a permutation of length $n$, to produce it by passing an ordered sequence through the stacks, do the following.  Define $I$ to be the input list, and initialize it to $1,\ldots, n $.
While $I\neq \emptyset$ do:
\be
\item If the top letter of $I,A$ or $B$ is next output letter, then output.
\item Else if $A=\emptyset$ then place next input letter on $A$.
\item Else apply Rules \ref{rules}. If two contradictory rules apply, then reject.
\ee

\noi
When $I=\emptyset$ then accept.
\end{alg}

For example, in Figure \ref{fig:stacks} we generate $52314$ by using rules 1.2 for step 
two, 3.1 for step three, and 2.2 for step four. The remaining steps simply offload the next tokens to be output.

\begin{figure}[ht!]
\bc
\bt{|c|c|c|}
\hline
\includegraphics[height= 35.0mm]{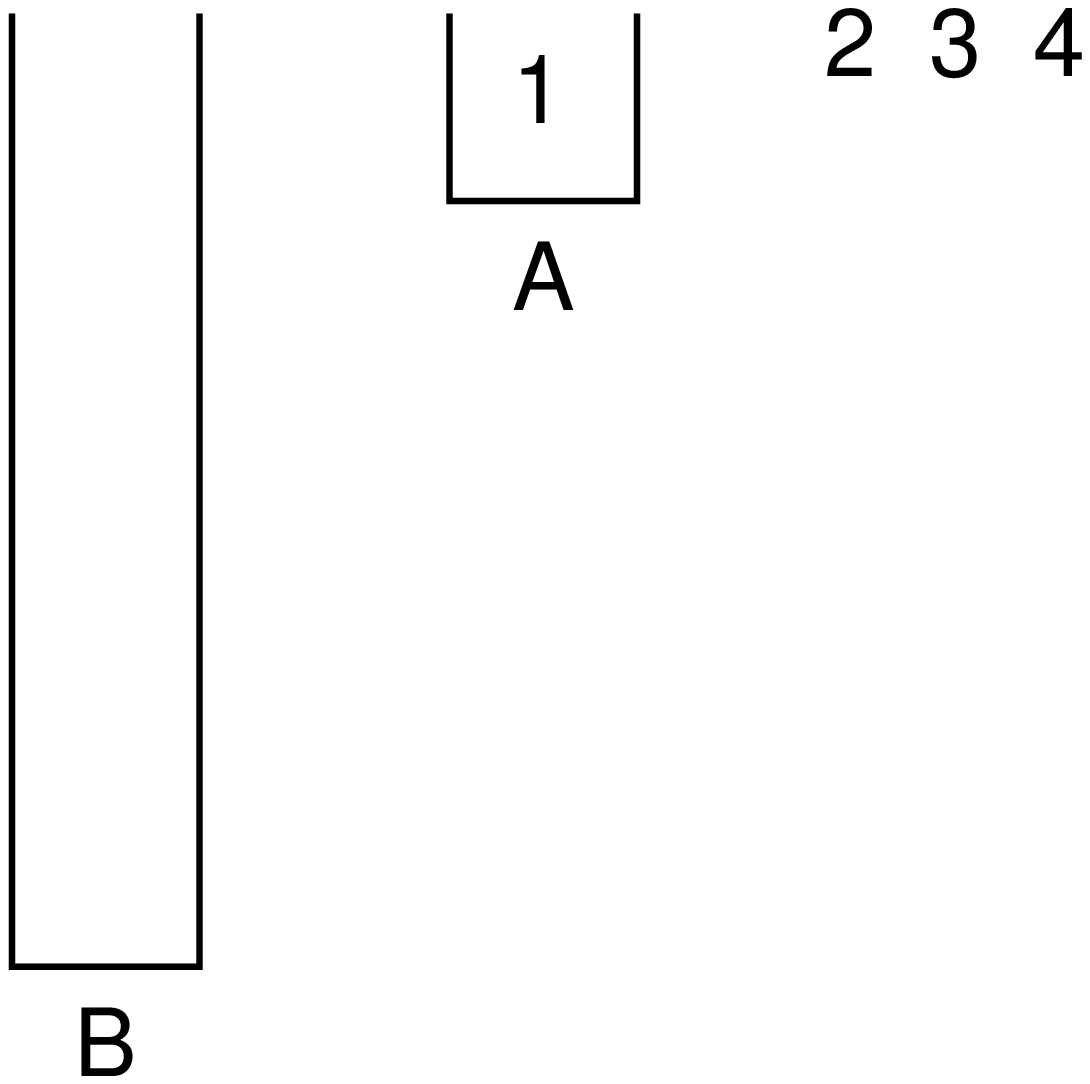} &
\includegraphics[height= 35.0mm]{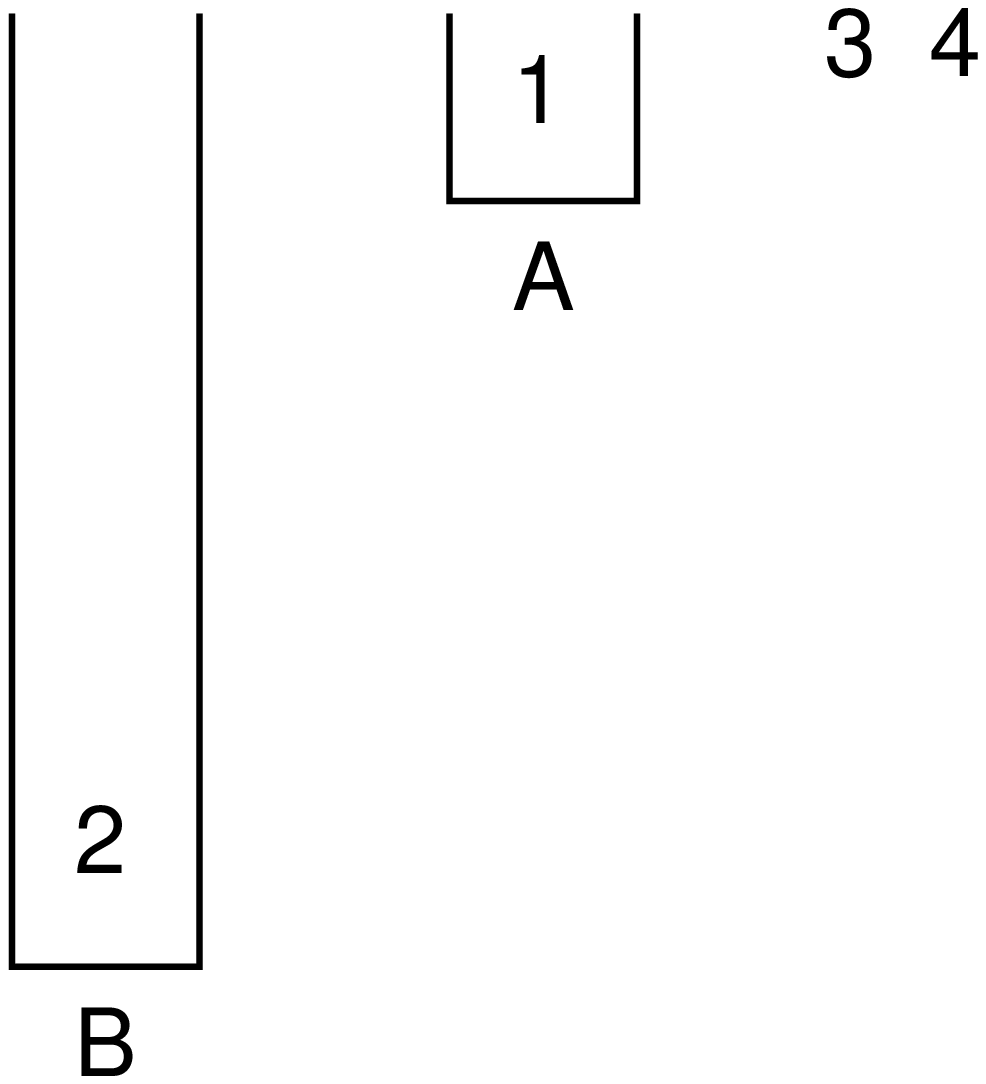} &
\includegraphics[height= 35.0mm]{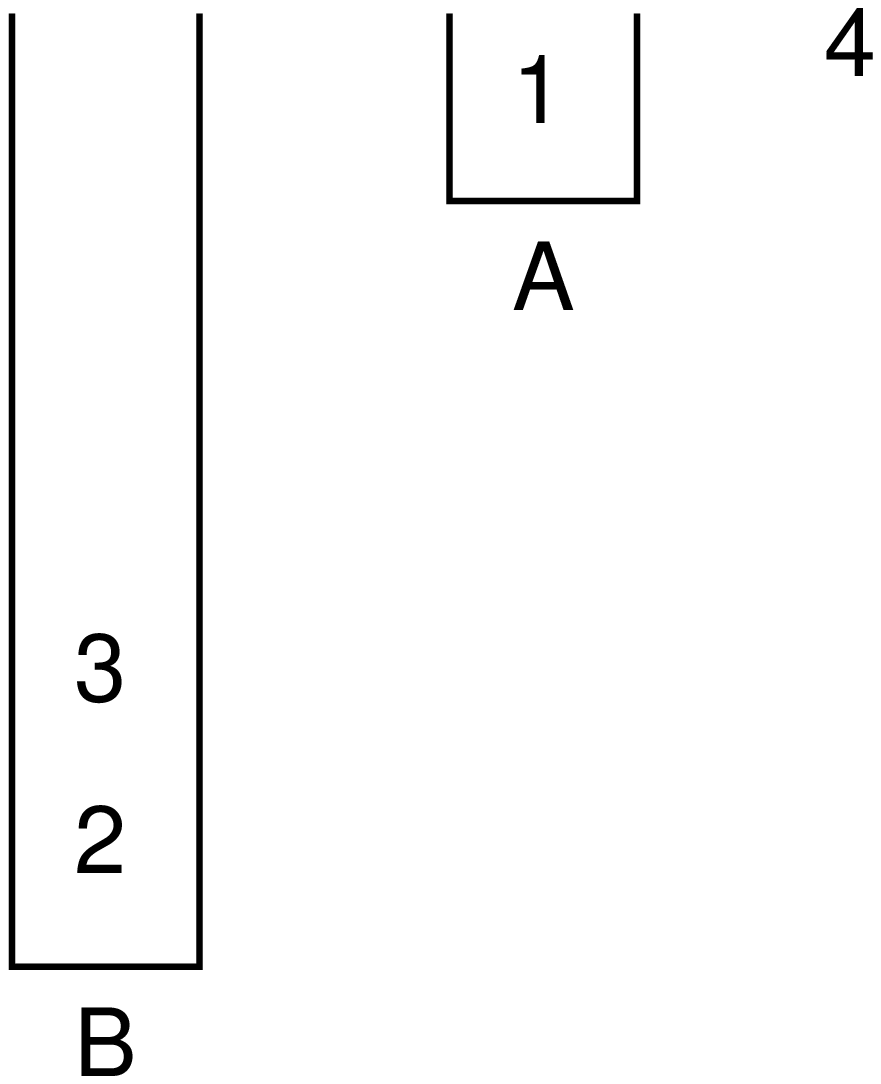} \\
\hline
 \texttt{Apply rule 3.1} &   \texttt{Apply rule 3.1} & \texttt{Apply rule 3.1}  \\
\hline
\et \ec
\caption{Generating the permutation $4132$.}\label{fig:stacks2}
\end{figure}

To generate the permutation $4132$ we keep $1$ on stack $A$ and pass $2$ to stack $B$, since the subinterval $41$ does not contain $3$ so is not right-contiguous modulo $2$. We show this in Figure  \ref{fig:stacks2}.

\begin{figure}[ht!]
\bc
\bt{|c|c|c|}
\hline
\includegraphics[height= 33.0mm]{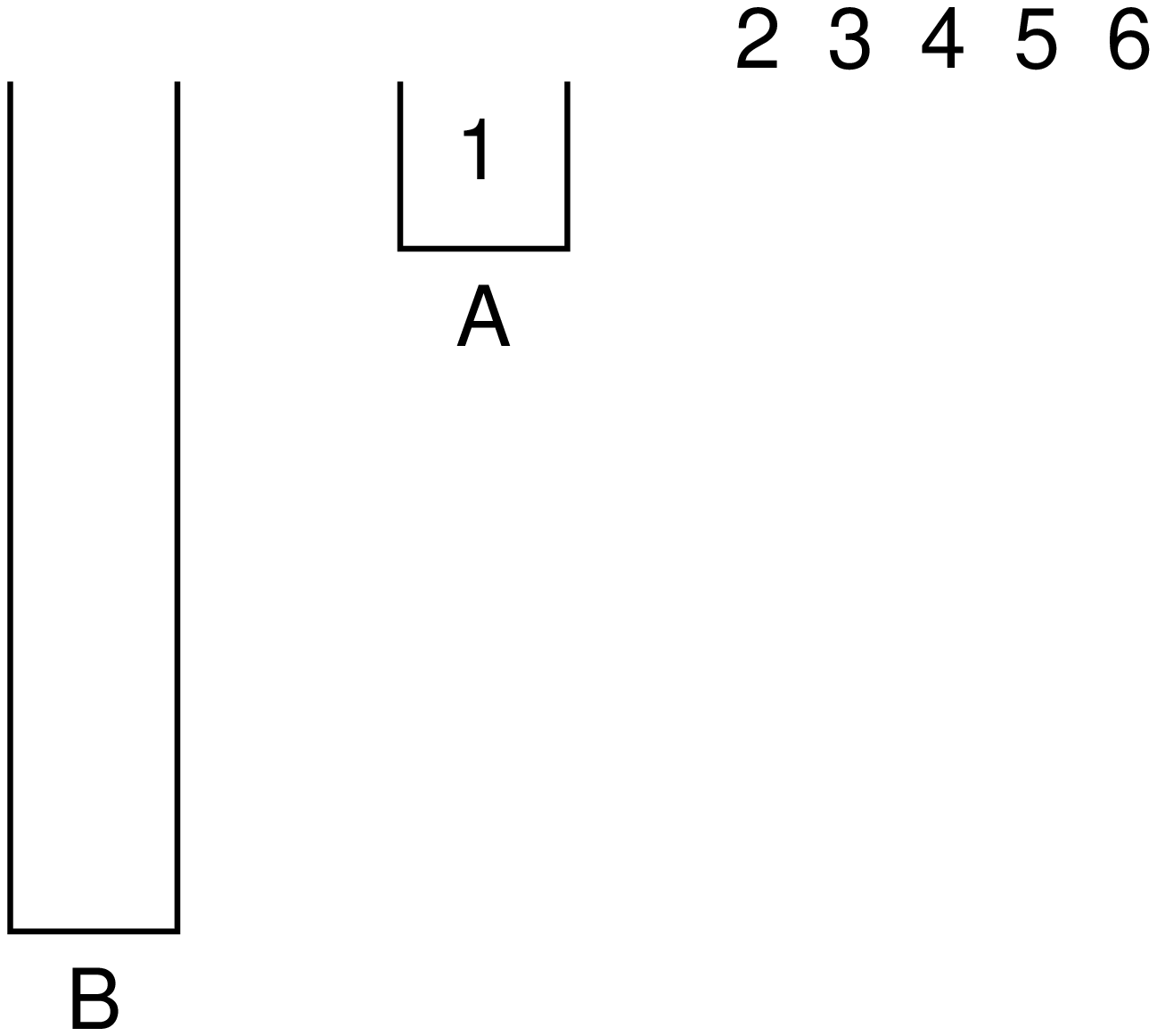} &
\includegraphics[height= 33.0mm]{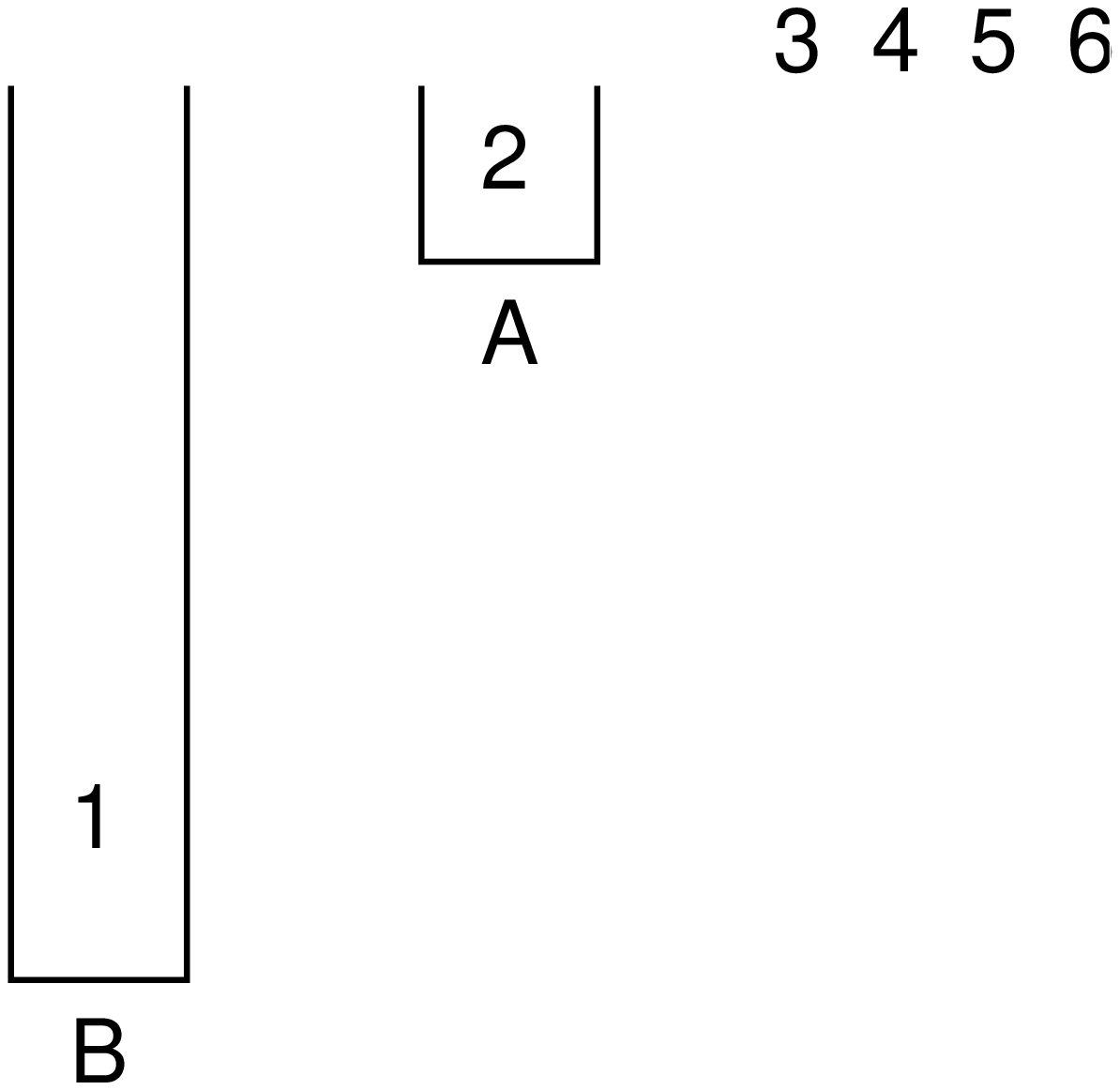} &
\includegraphics[height= 33.0mm]{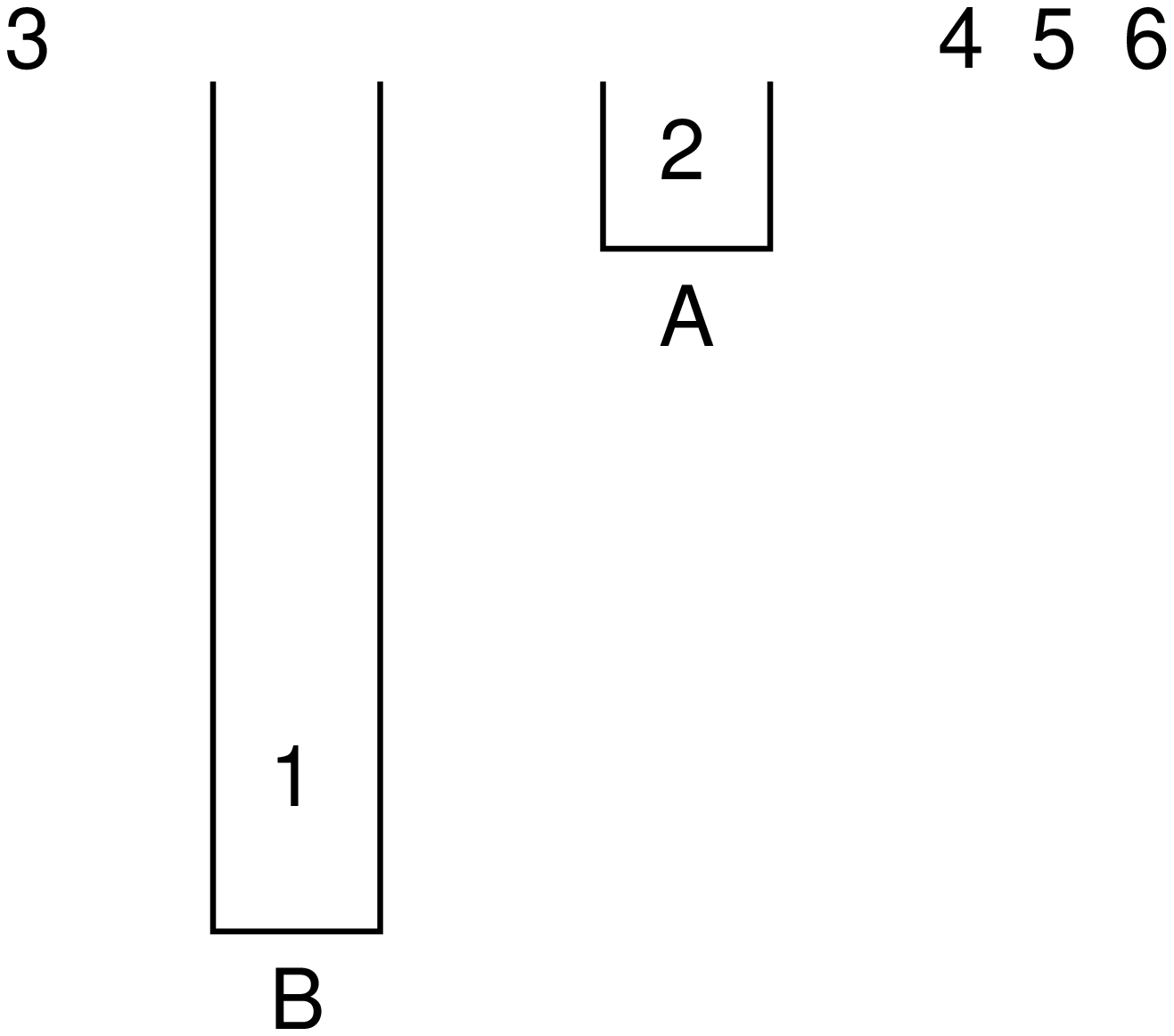} \\
\hline
 \texttt{Apply rule 3.1} & & \\
\hline
\includegraphics[height= 33.0mm]{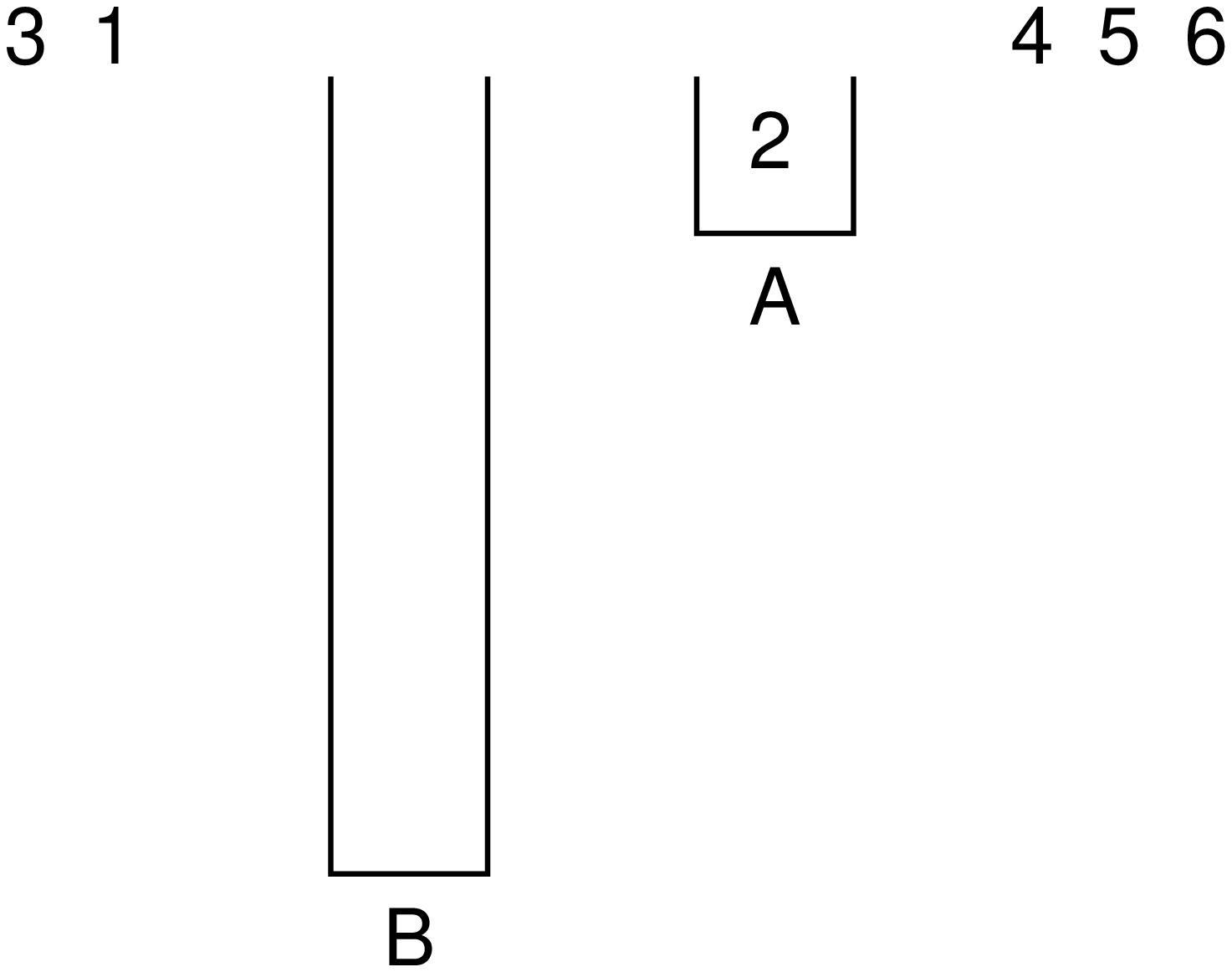} &
\includegraphics[height= 33.0mm]{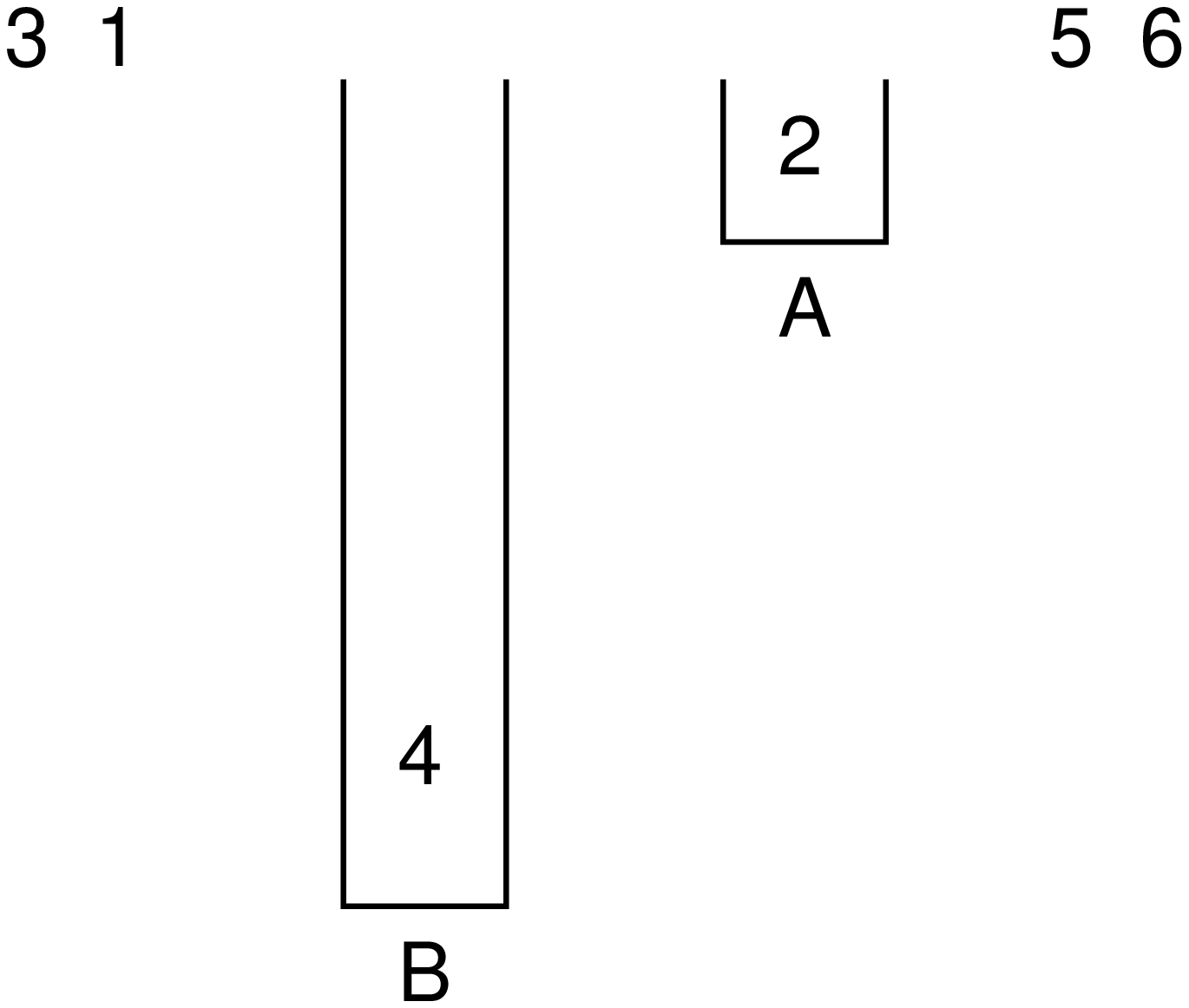} &
\includegraphics[height= 33.0mm]{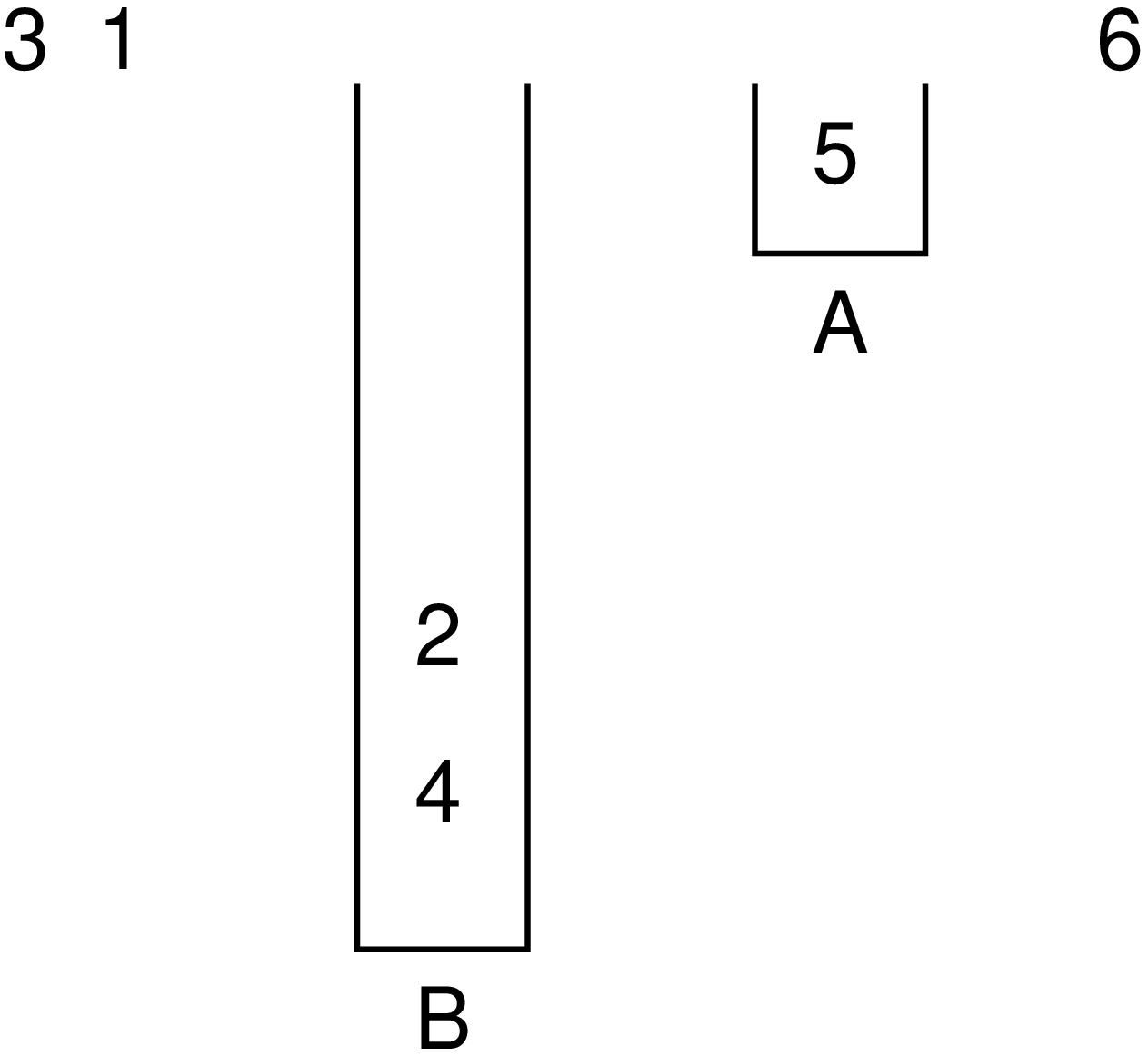} \\
\hline
 \texttt{Apply rule 3.1} & \texttt{Apply rule 2.2} & \\
\hline

\et \ec
\caption{Generating the permutation $316245$.}\label{fig:stacks3}
\end{figure}

As another example, consider the permutation $316245$. If you were to keep $1$ on 
stack $A$ and pass $2$ over to stack $B$ (contrary to rule 3) you would not succeed. 
The steps dictated by rules 1-3 are used in Figure \ref{fig:stacks3}. Compare the two applications of rule 3.1.

We will prove that if the permutation being generated avoids $\mathcal B$, then
 no two contradictory rules ever apply, and in this case that each rule preserves the well ordering of $B$. 
Therefore if a permutation is in $Av(\mathcal B)$ then it can be generated using this algorithm.

\begin{lem}\label{lem:wellordered}
If no two rules apply simultaneously then each move keeps the stack $B$ well ordered.
\end{lem}
\textit{Proof:}
 If rule 2.1 applies then there is a $b$ on $B$ that precedes $a$, but does not precede $x$
  or rule 2.2 would also apply. Similarly, if rule 2.2 applies then there is a $b$ on $B$ that 
  precedes $x$, but does not precede $a$ or rule 2.1 would also apply. In either case the 
  token that goes on $B$ precedes the token $b$. If neither rules 2.1 or 2.2 apply then $a$ 
  and $x$ both precede all tokens that may be on $B$, so we can place either on $B$, as directed by a rule 1.1, 1.2, 3.1 or 3.2, and 
  $B$ will remain  well ordered.
\hfill$\Box$

\begin{lem}\label{lem:cdlemma}
If a configuration where $c$ is the top entry of $B$ and $A$ is momentarily empty (before the next input letter fills it) is reached part way through producing an output permutation $\sigma$, where 
no two contradictory rules applied up to this point, then $\sigma$ has a sub-permutation which is order isomorphic to $4132$ or $4231$.
\end{lem}
\textit{Proof:}
Consider the point at which $c$ is placed on stack $B$. Some rule must have applied, with $c$ either next to be input with some  $d$ on stack $A$, or vice versa,  and some token $z$ in the input the next to be output in the permutation. The output is $z-d-c$ since we must get to the configuration with $c$ on 
$B$ and $A$ empty.

If rule 1 applied then there must be $p,q$ in the input before $z$ and output as $z-d-p<^q_c$. So taking 
$z=4,p=3$ and $c=1,d=2$ or $d=1,c=2$ we get $4d3c$.

If rule 2 applied to put $c$ on $B$ and $d$ on $A$ then there must be some $e$ on $B$ with the output of the form $z-c-e-d$, but $d$ must be output before $c$ to get the required configuration. So this is a contradiction.

If neither rules 1 or 2 apply, then since $d$ precedes $c$ in the output, rule 3 tell us to put $d$ on stack $A$ unless $z\ldots d$ is right-contiguous modulo $c$ and avoids $312$. Since $d$ is kept on $A$ while $c$ goes on $B$, we conclude that these conditions were not met.

If $z\ldots d$ fails to be right-contiguous because of some $p$ input after $c$ and $d$ and before $z$ and output after $d$ then we have output $z-d-p-c$ since $p$ must be output before the configuration with $c$ on $B$ and $A$ empty is reached. This is order isomorphic to $4d3c$ with $c=1,d=2$ or$c=2,d=1$. So we can assume that all tokens input between $c,d$ and $z$ precede $d$ in the output, in which case they must be in descending order in the output since $d$ occupies stack $A$ as they are input onto $B$.

If $z\ldots d$ contains a $312$ then since $d$ is the minimum entry in the subinterval and occurs last, the entries that make the $312$ are greater than $c$, so setting $c=1$ we have $4231$.

So the subinterval has no $312$ and the entries less than $z$ form a decreasing sequence, so the only other way it can fail the conditions is for there to be some $m>r>z$ in the input so $m$ is in $z\ldots d$ and $r$ is output after $d$. Since we reach the configuration where $c$ is on $B$ and $A$ is cleared, we know that $r$ is output before $c$, so we have $y-d-r-c$ which is $4d3c$ where $c=1,d=2$ or$c=2,d=1$. 
\hfill$\Box$

\begin{lem}\label{lem:contradict}
If a permutation avoids the list $\mathcal B$ then no two contradictory rules apply.
\end{lem}

We prove this by considering case-by-case when two contradictory rules apply at some instant, and in each case show that this occurs if the permutation being generated contains one or more of the permutations in $\mathcal B$. In Table \ref{table:Lem} we summarize these cases, and the reader may wish to refer to this as they read through the proof. Note that Rules 1.1 and 2.1 give the same instruction, so if both apply simultaneously they do not contradict each other. Similarly for Rules 1.2 and 2.2.

\textit{Proof:}
Suppose the permutation is part way generated as in Figure \ref{fig:2stacks}, and you
 reach a point where two contradictory rules apply for the first time. Since no two contradictory 
 rules applied before this,  by Lemma \ref{lem:wellordered} the stack $B$ is well ordered up to 
 this point.

\vspace{3mm}

\noi\textbf{Case: 1.1 and 1.2 apply.}
If 1.1 and 1.2 both apply,
then there exist $p,q$ in the input between $x$ and $y$ with either $p<q$ or $q<p$ so 
that $y-a-p-x$ or $y-x-a-p-q$ is output, which is one of $51234,52134,51243,52143$.

\vspace{3mm}

\noi\textbf{Case: 1.1 and 2.2 apply.}
If 1.1 and 2.2 both apply, then
 there exist $p,q$ in the input between $x$ and $y$ with either $p<q$ or $q<p$ so that $y-a-p-q$ is 
 output, and there exists $b$ on stack $B$ so that $b$ precedes $x$. Note that 1.2 does not apply 
 (covered by previous case) so $x$ comes after $p$, so we have $y-a-p<^q_x$ output.
If $b$ precedes $p$ then taking $(a,b)=(1,2)$ and $(p,q)=(3,4)$ we have one of $51234,52134,51243,52143$.
If $b$ is output after $p$ then taking $(a,b)=(1,2),x=3 $ and $p=4$ we have one of $51423,52413$.

Note that so far we have covered all basis permutations of length $5$.

\vspace{3mm}

\noi\textbf{Case: 1.2 and 2.1 apply.}
If 1.2 and 2.1 both apply, then there exist $p,q$ in the input between $x$ and $y$ with either $p<q$ or $q<p$ 
so that $y-x-p-q$ is output, and there exists $b$ on stack $B$ so that $b$ precedes $a$. Note that 1.1 
does not apply (covered by previous case) so $a$ comes after $p$, so we have $y-x-p<^q_a$ output.
If $b$ precedes $p$ then  taking $(x,b)=(1,2)$ and $(p,q)=(3,4)$ we have one of $51234,52134,51243,52143$.
If $b$ is output after $p$ then we have $y-x-p<^q_{b-a}$ as output. Go back to the instant that $a$ first goes on stack $A$ from the input.

If  $b$ was on stack $B$, then when $b$ was put there, there was a $c$ on $A$ and a $z$ in the input that was next to be output. So we have as output $z-c-y-x-p<^q_{b-a}$. Taking $z=4$, $(b,c)=(1,2)$, $y=7$, $x=5$, $p=6$ and $a=3$ we have either $4175623$ or $4275613$.

If $b$ was on stack $A$ when $a$ was input, if rule 2.2 applied then there must have been a token $d$ on 
$B$ that preceded $a$, and since no contradictory rules occur here then $b$ precedes $d$.
So we have as output $y-x-p<^q_{b-d-a}$. Taking $(b,d)=(1,2)$, $a=3$, $x=4$, $p=5$ and $y=6$
we get $645123$ and $645213$. 

If rule 1.2 applied then there must be two tokens $r,s$ between $a$ and the next output (which could be $y$ or another token earlier in the input) such that $a$ precedes $r,s$ in the output. But since $b$ precedes $a$ then $b$ would also precede $r,s$, so rule 1.1 would simultaneously apply, which is a contradiction.

\begin{figure}[ht!]
\bc
\includegraphics[height=4cm]{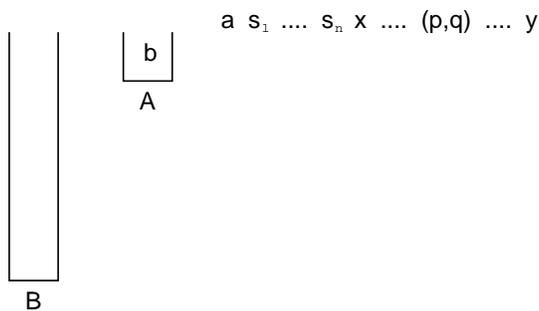}
\ec
\caption{The point that $a$ first goes on $A$, when $b$ is on stack $A$.}\label{fig:rule3}
\end{figure}

So if neither rules 1 or 2 applied, then rule 3 must apply. We have $b$ on $A$, $a$ next to be output, then possibly some $s_1,\ldots, s_n$ between $a$ and $x$, all of which will be output before $y$, then $x$, $p$ and $q$ in either order, then $y$ in the input. See  Figure \ref{fig:rule3}.
Since $b$ precedes $a$ in the output  it would stay on $A$ unless the subinterval from the next output letter $y$ all the way down to $b$ is  right-contiguous modulo $a$ and avoids $312$.
 Either there are no $s_i$ so $y$ is the next output letter, or if there are $s_i$ then they are all output before $y$, so the subinterval contains the segment from $t$ to $b$. But this is of the form $y-x-p-b$ where $x,p<y$ and $x<p$, so $y-x-p$ form a $312$ subsequence, so the subinterval does not avoid $312$.

Note that up to here  we have covered all basis permutations of length $5$ and the permutations 
$645123$, $645213$, $4175623$ and $4275613$.

\vspace{3mm}

If 2.1 and 2.2 both apply, then there exists $b$ on top of stack $B$ with 
$y-b<^{a}_{x}$ in the output. We consider two cases, when $a<b$ and $b<a$.

\vspace{3mm}

\noi\textbf{Case: 2.1 and 2.2 apply and $a$ precedes $b$.}
Suppose $a$ precedes $b$ in the input. Consider the point at which $a$ is on $A$ and $b$ is next to be input. If rule $1.1$ applies then we have some $p,q$ and next output $z$ (possibly $y$) and for output we have $z-a-p<^q_b$. But $b$ precedes $a$ so this case does not apply.

If rule $2.1$ applies then there exists $c$ on top of $B$ with output of the form $y-b<^{c-a}_x$. Remember that as input $a<b<x<y$. We may ask how $c$ came to be on stack $B$. If at some previous time $c$ is on top of $B$ and $A$ is momentarily empty (before the next input letter comes to fill it), then by Lemma \ref{lem:cdlemma} it must be that $4d3$ is output at the point that $c$ is on $B$ and $A$ is empty, then $y-b<^{c-a}_x$ is output, which is one of $y-b-c-a-x$, $y-b-c-x-a$ or $y-b-x-c-a$. These give 
$4d3ycax=4d37c56$, $4d3ycxa=4d37c65$ or $4d3ybxca=4d3867c5$ with $c,d=1,2$.

If neither rules $1.1$ or $2.1$ apply then rule 3 applies. If the next output letter is some $z<x$ then since $b$ precedes $a$ in the output, $b$ is placed on $A$ unless $z\ldots b$ is right-contiguous modulo $a$ and avoids $312$. But since $y$ is output between $z$ and $b$ and $x$ is output after, the subinterval is not right-contiguous, so $b$ would go on $A$, a contradiction.
   
  If the next output is greater than $x$ it must be $y$. Again the subinterval $y\ldots b$ is not right-contiguous since $x$ is output after $b$, so $b$ should be kept on $A$ and $a$ on $B$, a contradiction. 
   
\vspace{3mm}

\noi\textbf{Case: 2.1 and 2.2 apply and $b$ precedes $a$.}
Suppose that $b$ precedes $a$ in the input. Consider the instant that $a$ is to be input. If $b$ is on stack $A$, then some rule pushes $a$ onto $A$ and $b$ onto $B$. If rule $1.2$ is responsible, then there exist $p,q,z$ in the input with $z$ next to be output (possibly $z=y$) and $z>p>a$, $z>q>a$, and the output is $z-a-p<^q_b$. But $b$ precedes $a$ in the output so  this is a contradiction.

If rule $2.2$ is responsible then there exists $c$ on stack $B$ such that $b$ precedes $c$ and $a$ comes after $c$ in the output. Note that by assumption no two contradictory rules have applied before the point that $x$ is to be input. Thus for output we have $y-b<^{c-a}_x$ which is one of $5bc34$, $5bc43$ or $5b4c3$ with $b,c=1,2$.

If neither rules 1.2 nor 2.2 applied then rule 3 was responsible for pushing $a$ on $A$ and $b$ off $A$ to $B$. Now $b$ precedes $a$ in the output, so we must check whether the subinterval from the next token output to $b$ is  right-contiguous modulo $a$ and avoids $312$. Suppose the input sequence is 
$ap_1\ldots p_nx\ldots y$. If $y$ is next to be output then the subinterval $y\ldots b$ does not contain $x$, so fails to be right-contiguous. Similarly if some $p_i$ is the next to be output then $p_i\ldots b$ contains $y$ but not $x$, so fails to be right-contiguous.  So rule 3 would require that $b$ stay on stack $A$, which did not occur so this is a contradiction.

If $b$ is on stack $B$ and the token on stack $A$ is output before $a$ is input, then by Lemma \ref{lem:cdlemma} the output is of the form $4d3-y-b<^a_x$ for some $d$ with $b,d=1,2$, which is $4d37b56$ or $4d37b65$ .

Lastly,  $b$ is on stack $B$ and there is some
token $c$ on stack $A$ which cannot be output before $a$ is input, but is output by the time $x$ is to be input, so there is some $p$ next to be output with $a<p<x$. Since $c-b-a$ is the output order then rule $2.2$ applies to put $a$ on $A$ and $c$ on $B$ on top of the token $b$. So in full the  output is $p-c-y-b<^a_x$ where $p=4$, $a=3$, $x=5$, $y=6$ and $b,c=1,2$ so we have one of $4c6b35$ or $4c6b53$, and these are the last permutations in $\mathcal B$ to be accounted for.

\vspace{3mm}

Finally if none of rules $1$ or $2$ apply, then either rule $3.1$ or $3.2$ applies 
depending on whether $a$ precedes $x$ or not in the permutation. So we cannot 
have both $3.1$ and $3.2$ simultaneously.
\hfill$\Box$

We summarize the preceding proof in Table \ref{table:Lem}, which shows which  sub-permutations are forced
when two contradictory rules apply.
\begin{table}
\bc
\bt{|l|l|l|}
\hline
Rules 			&  Extra condition	 	& Sub-permutations responsible \\
\hline\hline
$1.1$ and $1.2$ 	& 		& $51234, 52134, 51243, 52143$\\
\hline
$1.1$ and $2.2$ 	& 		& $51234, 52134, 51243, 52143,51423,  52413$\\ 
\hline
$2.1$ and $1.2$ 	&		& $51234, 52134, 51243,  52143$,\\
			 	&		& $  645123, 645213, 4175623, 4275613$\\ 
\hline
$2.1$ and $2.2$ 	& $a<b$ 	& $4137256, 4237156, 4137265, 4237165,$\\
				& 		& $  41386725, 42386715$\\
\hline
				& $b<a$ 	& $51234, 52134,  51243, 52143, 51423, 52413$,\\ 
				&		& $416235, 426135, 416253, 426153,$\\ 
				&		& $4137256, 4237156,  4137265, 4237165 $\\
\hline

\et
\caption{Summarizing the proof of Lemma \ref{lem:contradict}.}\label{table:Lem}
\ec\end{table}

\begin{lem}\label{lem:avoidBimpliesAlg}
If a permutation avoids $\mathcal B$ then it can be generated by Algorithm \ref{alg}.
\end{lem}
\textit{Proof:}
If the algorithm halts while there is still input, it must be because of a contradiction between the rules $1-3$. By Lemma \ref{lem:contradict} this doesn't occur if the permutation avoids $\mathcal B$. If the algorithm does not halt until there is no more input, then since $B$ is well ordered by Lemma \ref{lem:wellordered} the permutation can be successfully generated.
\hfill$\Box$

Putting Lemmas \ref{lem:avoidBimpliesAlg} and \ref{lem:basis} together we get:
\begin{thm}\label{thm:main}
Let $\sigma$ be a permutation. The following are equivalent.
\be
\item $\sigma$ can be generated by a stack of depth two and an infinite stack.
\item $\sigma$ can be generated by Algorithm \ref{alg}.
\item $\sigma$ avoids the set of 20 permutations $\mathcal B$.
\ee
\end{thm}
\textit{Proof:} 
If  $\sigma$ avoids $\mathcal B$ then by Lemma  \ref{lem:avoidBimpliesAlg} it can be generated by  Algorithm \ref{alg}, so it can be generated by the stacks. If $\sigma$ contains a permutation from $\mathcal B$ then it cannot be generated by the stacks by Lemma  \ref{lem:basis}, so conversely if it can be generated by the stacks, it must avoid $\mathcal B$.
\hfill$\Box$

\section{Acknowledgments}

The author wishes to thank Steve Linton for suggesting the problem and conjecture, the reviewer for their detailed reading of the article and many helpful corrections and changes, as well as  Andrew Rechnitzer, Nik Ru\v{s}kuc, Vince Vatter, Steve Waton and Mike Zabrocki for fruitful discussions concerning this work. This work was in part supported by an EPSRC grant GR/S53503/01. I wish to dedicate this paper to my father Mr Brian Elder.


\begin{thebibliography}{10}

\bibitem{Atk}
M.Atkinson, M. Livesey and D. Tulley.
\newblock {\em Permutations generated by token passing in graphs}. Theor. Comp. Sci. 178 (1997), 103-118.

\bibitem{bona-stacks} Mikl\'os B\'ona. {\em Combinatorics of permutations}. Chapman \& Hall/CRC 2004. 


\bibitem{bona-book}
Mikl\'os B\'ona.
\newblock {\em A survey of stack-sorting disciplines}.
Electron. J. Combin. 2002-3 v.9 no.2.


\bibitem{me-vince}
Murray Elder and Vince Vatter.
\newblock {\em Problems and conjectures presented at the Third International Conference on Permutation Patterns, University of Florida, 2005.}
\texttt{arXiv.org/abs/math/0505504}


\bibitem{Knuth}
Donald Knuth.
\newblock {\em The art of computer programming: Sorting and searching}.
v.3 ed.2. Reading, Massachusetts: Addison-Wesley, 1998.


\bibitem{sal}
Steve Linton.
\newblock GAP program to study classes of permutations generated by by token passing networks via finite state automata. \texttt{http://www-groups.dcs.st-and.ac.uk/$\sim$sal}


\bibitem{Max}
Max Murphy.
\newblock {\em Restricted permutations, antichains, atomic classes, and stack sorting}. Doctoral Thesis, University of St Andrews, 2002.

\end{thebibliography}
\end{document}